%% file: m4-23.tex
\documentclass{gtart_h}

\input gtmonout

\volumenumber{4}
\volumename{Invariants of knots and 3-manifolds (Kyoto 2001)}
\volumeyear{2002}
\papernumber{23}
\pagenumbers{363}{376}
\received{27 June 2003}
\revised{31 March 2004}
\accepted{12 April 2004}
\published{2 May 2004}

\usepackage{amssymb} 
\usepackage{amsmath,amscd, epsf}

\newtheorem*{theorem}{Theorem}

\theoremstyle{definition}

\newcommand{\bD}{{\mathbb D}}

\newcommand{\bQ}{{\mathbb Q}}

\newcommand{\bS}{{\mathbb S}}

\newcommand{\bZ}{{\mathbb Z}}

\begin{document}
\title{Some computational results on mod 2 finite-type\\invariants 
of knots and string links}                    
\shorttitle{Some computational results on mod 2 finite-type invariants}
\authors{Ted Stanford}                  
\address{New Mexico State University, Dept of Mathematical Sciences,
PO Box 30001\\Department 3MB, Las Cruces, New Mexico 88003-8001, USA}

\email{stanford@nmsu.edu}                     
\begin{abstract}   
We publish a table of primitive finite-type invariants of order
less than or equal to six, for knots of ten or fewer crossings.
We note certain mod-2 congruences, one of which leads to 
a chirality criterion in the Alexander polynomial.  We state
a computational result on mod-2 finite-type invariants of
2-strand string links.
\end{abstract}


\maketitle 

\section{Introduction}

In \cite{Vassiliev:CohomologyKnot1990}, 
Vassiliev described a new way to obtain invariants of knots in
$\bS^3$.  His paper contains the outline of an algorithm for computing
his invariants.  Gusarov \cite{Gusarov:nEquivalence1994} 
obtained the same set of invariants
independently and by different methods.  These invariants of
Vassiliev and Gusarov are now often referred to as finite-type
invariants.

At the end of this paper there are tables of primitive invariants of
order $\le 6$ for knots of $\le 10$ crossings.  These invariants were
computed using an implementation of Vassiliev's algorithm, which is
described in \cite{Stanford:ComputingVassiliev1997} .  In order to
create tables such as these, a basis must be chosen for the invariants
of order $\le 6$.  In Section~\ref{S:basis} we will discuss the choice
of such a basis, and make some related observations.  In
Section~\ref{S:string_links} we note that the algorithm for computing
knot invariants extends easily to the computation of finite-type
invariants of string links.  We describe a computation which shows
that there is a mod-2 weight system of order $5$ (first noted by
Kneissler and Dogolazky) for $2$-strand string links which does not
``integrate'' to a mod-2 finite-type invariant of order~$5$.  In
Section~\ref{S:table} we present the two matrices for translating our
numbers into finite-type invariants obtained from the derivatives of
knot polynomials, following the notation of
Kanenobu~\cite{Kanenobu:Order6_1998}.

\section{Choosing a basis}\label{S:basis}

For a general reference on finite-type invariants, see
Bar-Natan~\cite{Bar-Natan:OnVassiliev1995} 
or Birman~\cite{Birman:NewView1993}.

Any $\bQ$-valued knot invariant $v$ may be extended to singular knots
in a unique way by the usual formula:
\begin{equation}
v(K_\times) = v(K_+)-v(K_-)\
\label{ccformula}
\end{equation}
If there exists a positive integer $n$ such that $v(K) = 0$
for any knot $K$ with more than $n$ double points, then $v$
is said to be a finite-type invariant.  The least such $n$
is called the order or type or degree of $v$.  Let $V_n^*$
be the $\bQ$ vector space of knot invariants of order $\le
n$.  $V_n^*$ is finite dimensional because there are a
finite number of chord diagrams with $n$ chords or less.
Let $W^*_n \subset V_n^*$ be the subspace of invariants
which are additive under the connected sum of knots.  (That
is, $w(K_1\#K_2) = w(K_1) + w(K_2)$ for all $w \in W^*_n$
and for all knots $K_1,K_2$.)  Additive finite-type
invariants are often referred to as primitive invariants
because they are the primitive elements of the graded Hopf
algebra $\cup_{i=0}^\infty V_n^*$.  This means that all
finite-type invariants are linear combinations of products
of the primitive ones.  Therefore, in making tables of
invariants, it suffices to list only basis elements for
$W^*_n$.

There does not seem to be a canonical way to choose a basis for $W^*_n$.
We list some of the desirable properties that such a basis 
$B_n = \{b_1, b_2, \dots b_{{\rm dim}(W^*_n)}\}$ might have:

\begin{enumerate}
\item
$B_n$ should consist of $\bZ$-valued invariants.

\item\label{integer_basis}
$B_n$ should be a basis over $\bZ$ for the primitive integer-valued
invariants of order $\le n$. 

\item
For each $i \le n$, the set $\{b_1, b_2, \dots b_{{\rm dim}(W^*_i)}\}$ should
be a basis for $W^*_i$.  

\item
Knots of small crossing number should have
small values on the basis invariants.  
I first did the computations presented here
in 1992, and at that time an essentially random
basis was obtained from the computer program that 
solved the T4T relations.
A table of invariants using
that basis was made available electronically, 
though it was never published.  That table
contained many four-digit numbers.
Recently I have been able 
by ad-hoc methods to change the basis to give
the values shown below, where the largest absolute
value occuring is 39.

\item\label{even-odd}
If $w$ is an even-order basis element, then
$w(m(K)) = w(K)$ for any knot $K$, where $m(K)$
denotes the mirror image of $K$.  
If $w$ is an odd-order
basis element, then $w(m(K)) = -w(K)$ for any
knot $K$.  This is always possible over $\bQ$
(in fact over any ring where $2$ is invertible, as noted
by Vassiliev~\cite{Vassiliev:CohomologyKnot1990}),
using the identity 
$w(K) = \frac 12 (w(K)+w(m(K))) + \frac 12 (w(K)-w(m(K)))$.
The first term on the right will always have even order, and
the second will always have odd order, so $w$ may be replaced
by one of the two terms (modulo invariants of lower order), 
depending on whether the order of $w$ 
itelf is even or odd.

\end{enumerate}

Even though each of the above conditions 
can be satisfied individually,
it is not possible to satisfy
them all at once.
The basis invariants below
satisfy all the conditions
except~\ref{integer_basis}
and they {\it almost} 
satisfy~\ref{integer_basis}.
The vectors in $\bZ^{12}$ that actually occur as the values
for specific knots form a sublattice of $\bZ^{12}$ of
index 16, reflecting four inevitable mod-2 congruences
imposed by Condition~\ref{even-odd}.
The basis is chosen so that 
that $v_3 \equiv v_{4a}$ modulo~2, and likewise
$v_{5a} \equiv v_{6a}$, $v_{5b} \equiv v_{6b}$, and
$v_{5c} \equiv v_{6c}$.  

We note that 
$v_{4a} = \frac12 (3a_2-a_2^2) + a_4$,
where $\sum a_i x^i$ (with $i$ even)
denotes the Conway polynomial
of a knot.  If $v_3(K)$ is odd, then $v_3(K) \ne 0$,
and therefore $K$ is chiral.
Hence if $\frac12 (3a_2-a_2^2) + a_4$
is odd (or, to make it slightly more simple,
if $\frac12 (a_2+a_2^2) + a_4$ is odd),
then $K$ is chiral.  It is
interesting to find this chirality criterion contained
in the Conway polynomial.  
Stoimenow~\cite{StoimenowSquareNumbers2000p}
has recently studied
the chirality information in the determinant of a knot,
which is the integer obtained by evaluating the Conway
polynomial at $x = 2\sqrt{-1}$.
The chirality criterion from $v_{4a}$ modulo~2 
is independent of 
determinant.  More precisely, given
a number $d$ which occurs as the determinant of some knot,
there exist knots $K,K^\prime$, each with determinant $d$,
such that $v_{4a}(K) \equiv 0$
and $v_{4a}(K^\prime) \equiv 1$, modulo 2.  To see this, note
that replacing $a_4$ by $a_4+1$ and $a_2$ by $a_2+4$ does
not change the determinant of a knot, but it does change
$v_{4a}$ modulo 2. Such a change can be made because any
even polynomial with constant term
equal to $1$ occurs
as the Conway polynomial of some knot.  

\section{Mod-$2$ invariants of $2$-strand string links}\label{S:string_links}

Let $\bD^2$ be the two-dimensional disk, and let $p_1, p_2 \dots p_k \in \bD^2$
be $k$ distinct points.
A $k$-strand string link is a $k$-tuple of 
disjoint tame
curves in $\bD^2 \times [0,1]$ such that the endpoints of 
the $i$th curve
are $(p_i,0)$ and $(p_i,1)$ for all $1 \le i \le k$.
The components of a string links are thus ordered,
and each component has an unambiguous orientation.  A string link may
be given by a planar diagram, just as in the case of knots, and the
three usual Reidemeister moves suffice to generate equivalence.  We
may also consider the larger set of singular string links, which are
allowed to contain a finite number of double point singularities, just
as in the case of knots. (Two extra Reidemeister moves are needed here,
see~\cite{Stanford:ComputingVassiliev1997}.)
Applying Relation~\ref{ccformula}, we may
define a string link invariant $v$ to be of finite-type if there exists
a positive integer $n$ such that $v$ vanishes on string links with
more than $n$ singularities.  The least such integer $n$ is called
the order or type or degree of~$v$.  

Rather than deal with finite-type invariants directly, it is
often convenient to consider the abelian group generated by
all singular string links, subject to Relation
\ref{ccformula} and to the relation that $v(L)=0$ if $L$ has
more than $n$ singularities.  (Both of these are of course
infinite families of relations.)  For $k$-component string
links, we denote this group by $V_n(k)$.  The set of
finite-type invariants of order $\le n$ taking values in an
abelian group $G$ is then identified with ${\rm Hom}(V_n(k),
G)$ in the obvious way.

Although $V_n(k)$ is defined by an infinite presentation, it is
well-known, and easily seen, that $V_n(k)$ is a finitely-generated
abelian group.  The following is a version
of the well-known
``fundamental theorem of Vassiliev invariants'' (see Bar-Natan and
Stoimenow~\cite{Bar-NatanStoimenow:FundamentalTheorem1997}), and is easy
to prove using a the same methods as in
\cite{Stanford:ComputingVassiliev1997}:

\begin{theorem}
\label{thm:fundamental}
A presentation of $V_n(k)$ is given by any
set $S$ which contains exactly one string link
for each chord diagram with $\le n$ chords,
subject to the topological
4-term and 1-term relations, exactly one such relation from each
configuration class of order $\le n$.
\end{theorem}

As with singular knots, to every singular string link there
corresponds a chord diagram which records the combinatorial
information of the order of occurence of the double points in the
string link.  To every  T4T or T1T relation there corresponds
a (combinatorial) 4T or 1T relation, obtained by replacing each singular
string link with its associated chord diagram. Two T4T or
T1T relations are said to have the same configuration class if their
associated 4T or 1T relations are the same.

In order to make sense of the above Theorem,
it is necessary to
understand how an arbitrary T4T or T1T relation can be considered a 
relation among the elements of $S$, since these may be chosen
in a completely different way from the relations.  In the case
of those string links $L \in S$ which have $n$ singularities,
and the relations of order $n$,
there is no problem.  If $L$ and $L^\prime$ have the same chord
diagram with $n$~chords, then $[L] = [L^\prime] \in V_n$.  Hence
the T4T and T1T relations among singular string links become
4T and 1T relations among chord diagrams.  (These are the
``top row'' relations of Vassiliev~\cite{Vassiliev:CohomologyKnot1990}.)

Now suppose we have a T4T or a T1T relation of order $k<n$
Suppose
$L$ is string link in the relation.  Then there exists an 
$L^\prime \in S$
such that $L$ and $L^\prime$ share the same chord diagram.  
It is then possible
to make crossing changes to $L$ until it is equivalent to $L^\prime$. 
Thus our relation of order $k$ becomes a relation among 
elements of $S$ of order $k$, {\it plus} a sum of singular string links of
order $>k$.  
Each of these higher-order string links may in turn
be written as a sum elements of $S$ plus higher-order singular string
links.  Inductively, we see that each T4T or T1T relation becomes
a relation among the elements of $S$.  But there is no reason
that the relations should be homogeneous with respect to the
degree of the elements of $S$.

(Part of the content of
the Theorem
is that it does not matter what sequence
of crossing changes you choose to transform
a singular string link $L$ to 
$L^\prime \in S$.  More specifically,
different choices of crossing change sequences
will produce relations which differ by higher-order
T4T relations.)

Let $A_n = A_n(k)$ be the abelian group generated by all
chord diagrams of order~$n$, subject to the 
4T and 1T relations.  We see that there is 
a homomorphism $\phi:A_n \to V_n$, where if $D$
is a chord diagram of order~$n$ then $\phi(D)$
is any singular string link of order~$n$ whose
chord diagram is~$D$.
The combinatorial 4T and 1T relations are easier to work with
than their topological counterparts, and it would be nice if
the map $\phi$ were always injective, indeed this injectivity
question goes back to Vassiliev~\cite{Vassiliev:CohomologyKnot1990}.
The Kontsevich integral gives an almost satisfactory answer to this
question.  It works as well for string links in general as
it does for knots, and the result may
be stated as follows:

\begin{theorem}[Kontsevich]\label{thm:Kontsevich}
Let $\phi: A_n(k) \to V_n(k)$ be as above.
Then the kernel of $\phi$ is finite.
\end{theorem}

The existence of torsion in $A_n(1)$ (the case of knots)
is unknown.  For the case of two-strand string links,
the element $\Delta w$, shown below,
was found by~Dogolazky
and~Kneissler (see~\cite{Dogolazky:Schlingeldiagramme1998})
to have order~$2$ in $A_5(2)$:

\bigskip
\epsfxsize=2in
\centerline{$\Delta w  = $ \epsffile {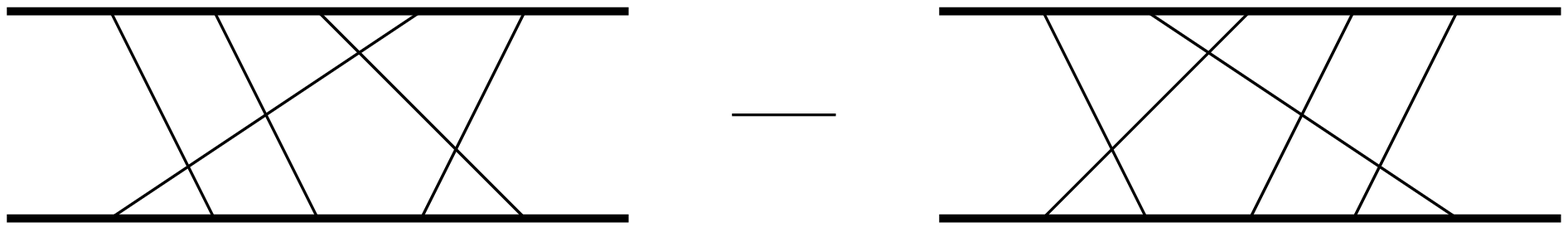}}
\bigskip

Using modified versions of the same computer programs
I used to make the tables at the end of this paper,
I have found that $\Delta w$
is in the kernel of $\phi$.
In fact, it is killed by T4T relations of order~4.
Thus there is no analog of the Kontsevich integral
over $\bZ/2\bZ$, at least for string links.

Here is a brief description of the computer programs
used to obtain this result.  There are two main programs.
The first program takes an abitrary singular string link
$L$ and makes crossing changes so that it is equivalent
to a standard $L^\prime$ with the same chord diagram.
In fact, the set $S$ of the standard links $L^\prime$ is not
kept explicitly.  Rather, there is a process of moves
and crossing changes to $L$ which results in a
``canonical'' $L^\prime$.  The process is a slight
modification of the algorithm described 
in~\cite{Stanford:ComputingVassiliev1997}.
First, a spanning tree is chosen for $L$,
where $L$ is viewed here as a spatial graph
with four-valent rigid vertices.  (The endpoints
of the string link may also be treated as edges
adjacent to a rigid vertex which is the boundary
of $\bD^2 \times [0,1]$.) Also,
a cyclic orientation is chosen for the
edges at each vertex, compatible with the
dihedral orientation required by the rigidness
of the vertex.  The choices of spanning tree and 
cyclic orientation are done by a simple but
arbitrary algorithm, whose only important property
is that it makes the same choices for any two
singular string links with the same chord diagram.

After the spanning tree and orientation are chosen,
Reidemeister moves are performed until the spanning
tree has no crossings on it, and until all the cyclic
orderings on all the vertices match those chosen
above.  Then the remaining edges of the spatial
graph $L$ are layered.  Each time a crossing change
is made, the new singular string link (with one more
singularity than $L$) is inductively processed by the
same algorithm until the number of singularities 
exceeds the order of the invariant to be computed.

The second program used in these computations is a generator
of T4T~relations.  A list of the configuration classes of
4T~relations is generated, and then the program realizes
each one as a T4T~relation.  As noted above, it does not
matter which particular T4T~relation is chosen, since any
two T4T~relations of the same configuration class will
differ by T4T~relations among string links with more
singularities.

After the T4T~relations are generated, they are fed into
the first program, and turned into linear combinations 
among the singular string links in the chosen set~$S$
(which exists only implicitly, as noted above).  The
result is a list of linear equations with the variables
indexed by the chord diagrams with $\le n$ chords.
The equations are solved, and a basis chosen.  
Thereafter, in order to compute the invariants
of a give string link $L$, only the first program
is necessary.

Unfortunately, there doesn't seem to be an easy way
to extract from the data files an understandable
linear combination of T4T relations which adds up
to $\Delta w$.  It would be nice to have an
understanding of how exactly the $2$-torsion is
killed.

\section{Knot tables}\label{S:table}

First we give matrices to translate our invariants into finite-type
the invariants obtained from standard knot polynomials, following
Kanenobu~\cite{Kanenobu:Order6_1998}.

The HOMFLYPT polynomial of a knot $K$ is written 
$$P(K; t,z) = \sum_{i=0}^NP_{2i}(K;t)z^{2i}$$
where 
$P_{2i}(K,t) \in \bZ [t^{\pm t}]$, and is 
determined by the skein relation
$$t^{-1}P(L_+; t,z)-tP(L_-; t,z) = zP(L_0; t,z)$$
The Jones polynomial $V(L;t)$ is given by
$$V(L;t) = P(L;t,t^{1/2}-t^{-1/2})$$
The Conway polynomial is given by $$\Delta_K(z) = P(K;1,z)$$
and is written
$$\Delta_K(z) = \sum_{i=0}^N a_{2i}(K)z^{2i}$$
The Kauffman polynomial of a knot $K$ is written 
$$F(K;a,z) = \sum_{i=0}^NF_i(K;a)z^i$$ 
where 
$F_i(K;a) \in \bZ [a^{\pm 1}]$, and is determined by
the skein relation 
$$aP(L;a,z) + a^{-1}P(L_-; a,z) = z(F(L_0;a,z) + a^{-2\nu}F(L_\infty; a,z))$$
The notation $P_{2i}^{(n)}$ denotes the knot invariant
obtained by evaluating the  $n$th derivative of the polynomial
$P_{2i}{K;t}$ at $t=1$, and similarly for the polynomials $V$ and $F_i$.
$$
\left[v_2, v_3, v_{4a}, v_{4b}, v_2^2, v_{5a}, v_{5b}, v_{5c}, v_2v_3, 
v_{6a}, v_{6b}, v_{6c},  v_{6d}, v_{6e}, v_2^3, v_3^2, v_2v_{4a}, v_2v_{4b}\right] M_1$$
$$=\left[a_2, a_4, \frac{P_0^{(3)}} {24},\frac{P_0^{(4)}} {24}, a_2^3, a_2a_4, 
\frac {a_2P_0^{(4)}} {24},
(\frac {P_0^{(3)}} {24})^2, \frac {V^{(5)}}{5!}, \frac {V^{(6)}}{6!}\right]
$$

where
\tiny
$$M_1=
\left[
\begin{array}{rrrrrrrrrr}
1 &  {\displaystyle -\frac32}  & 1 & -93 & 0 & 0 & 0 & 0 & 
128 & {\displaystyle -\frac {5327}{2}}  \\ [2ex]
0 & 0 & 2 & -12 & 0 & 0 & 0 & 0 & 275 & {\displaystyle -\frac {1345}{2}}  \\ [2ex]
0 & 1 & 0 & 24 & 0 & 0 & 0 & 0 & -30 & {\displaystyle \frac {1177}{2}}  \\ [2ex]
0 & 0 & 0 & -16 & 0 & 0 & 0 & 0 & 24 & -538 \\
0 & {\displaystyle \frac {1}{2}}  & 0 & 8 & 0 & {\displaystyle 
-\frac {3}{2}}  & -93 & 1 & -9 & 201 \\ [2ex]
0 & 0 & 0 & 0 & 0 & 0 & 0 & 0 & -54 & 135 \\
0 & 0 & 0 & 0 & 0 & 0 & 0 & 0 & 138 & -345 \\
0 & 0 & 0 & 0 & 0 & 0 & 0 & 0 & 30 & -75 \\
0 & 0 & 0 & 0 & 0 & 0 & -12 & 4 & -18 & 45 \\
0 & 0 & 0 & 0 & 0 & 0 & 0 & 0 & 0 & -81 \\
0 & 0 & 0 & 0 & 0 & 0 & 0 & 0 & 0 & 27 \\
0 & 0 & 0 & 0 & 0 & 0 & 0 & 0 & 0 & 21 \\
0 & 0 & 0 & 0 & 0 & 0 & 0 & 0 & 0 & -270 \\
0 & 0 & 0 & 0 & 0 & 0 & 0 & 0 & 0 & 156 \\
0 & 0 & 0 & 0 & 1 & {\displaystyle \frac {1}{2}}  & 8 & 0 & 0 & 
{-\displaystyle \frac {9}{2}}  \\ [2ex]
0 & 0 & 0 & 0 & 0 & 0 & 0 & 4 & 0 & 18 \\
0 & 0 & 0 & 0 & 0 & 1 & 24 & 0 & 0 & -45 \\
0 & 0 & 0 & 0 & 0 & 0 & -16 & 0 & 0 & 36 \\
\end{array}
\right]
$$\normalsize
and
$$
\left[v_2, v_3, v_{4a}, v_{4b}, v_2^2, v_{5a}, v_{5b}, v_{5c}, v_2v_3, 
v_{6a}, v_{6b}, v_{6c},  v_{6d}, v_{6e}, v_2^3, v_3^2, v_2v_{4a}, v_2v_{4b}\right] M_2 
$$
$$
=\left[\frac {P_0^{(6)}}{6!},\frac {P_2^{(4)}}{4!}, \frac {P_4^{(2)}}2, a_6, \frac {F_0^{(6)}}{6!},
\frac {F_1^{(5)}}{5!}, \frac {F_2^{(4)}}{4!}, \frac {F_3^{(3)}}{3!}, \frac {F_4^{(2)}}2, F_5^{(1)}
\right]
$$
where
\tiny
$$M_2=
\left[
\begin{array}{rrrrrrrrrr}
{\displaystyle -\frac {14065}{3}}  & 2307 & -348 & 
{\displaystyle \frac {40}{3}}  & {\displaystyle \frac {14065}{3}
}  & 9697 & {\displaystyle \frac {12413}{2}}  & {\displaystyle 
\frac {4478}{3}}  & 120 & 4 \\ [2ex]
-1075 & {\displaystyle \frac {533}{2}}  & -9 & 0 & -1075 & -1524
 & {\displaystyle -\frac {1133}{2}}  & -56 & -1 & 0 \\ [2ex]
1044 & {\displaystyle -\frac {1017}{2}}  & 72 & -2 & -1044 & 
-2136 & {\displaystyle -\frac {2719}{2}}  & -335 & -20 & 4 \\
 [2ex]
-888 & 389 & -48 & 1 & 888 & 1910 & 1375 & 425 & 56 & 2 \\
380 & -205 & 34 & {\displaystyle -\frac {3}{2}}  & -380 & -706 & 
{\displaystyle -\frac {917}{2}}  & -116 & -8 & 0 \\ [2ex]
240 & -72 & 3 & 0 & 240 & 320 & 84 & 0 & 0 & 0 \\
-560 & 144 & -5 & 0 & -560 & -800 & -288 & -20 & 1 & 0 \\
-160 & 60 & -3 & 0 & -160 & -192 & -12 & 16 & 1 & 0 \\
80 & -24 & 1 & 0 & 80 & 96 & 36 & 4 & 0 & 0 \\
-160 & 96 & -18 & 1 & 160 & 320 & 168 & 24 & 0 & 0 \\
32 & 0 & -6 & 1 & -32 & -80 & -96 & -68 & -18 & 0 \\
64 & -56 & 14 & -1 & -64 & -96 & 8 & 24 & -2 & -2 \\
-512 & 288 & -48 & 2 & 512 & 1056 & 608 & 104 & 0 & 0 \\
256 & -112 & 12 & 0 & -256 & -576 & -432 & -136 & -20 & -2 \\
{\displaystyle -\frac {32}{3}}  & 8 & -2 & {\displaystyle \frac {
1}{6}}  & {\displaystyle \frac {32}{3}}  & 16 & 8 & 
{\displaystyle \frac {4}{3}}  & 0 & 0 \\ [2ex]
32 & -16 & 2 & 0 & -32 & -64 & -48 & -16 & -2 & 0 \\
-96 & 64 & -14 & 1 & 96 & 176 & 88 & 4 & -4 & 0 \\
64 & -32 & 4 & 0 & -64 & -128 & -96 & -32 & -4 & 0
\end{array}
 \right] 
$$\normalsize
The notation used in the tables is as follows.  As above, 
$v_2$ is an invariant of order $2$, $v_{4a}$ and $v_{4b}$ are 
invariants of order $4$, and so forth.  Odd-order invariants
in the tables change sign under mirror image, and even-order
invariants are unchanged under mirror image.

The numbering of the
knots follows Rolfsen~\cite{Rolfsen:KnotsLinks1976}, 
up to mirror images.  The invariants were computed
directly from the braid words listed, which were
obtained from Jones~\cite{Jones:HeckeAlgebra1987}.
The knots 10.167 and 10.170 are switched with respect
to~\cite{Jones:HeckeAlgebra1987}, 
but not with respect to~\cite{Rolfsen:KnotsLinks1976}.
In each braid word, 
$a=\sigma_1$ (a positive crossing between the first and second
braid strands), $A=\sigma_1^{-1}$, $b=\sigma_2$, $B=\sigma_2^{-1}$,
and so forth.

\Addresses\recd
\newpage

\section*{Tables}
\addcontentsline{toc}{section}{Tables}

\def\newstrut{\vrule width 0pt
height 6pt depth 3pt}
\tiny
\cl{
\begin{tabular}{|l|l|r|r|rr|rrr|rrrrr|}
\hline\newstrut
Knot & Braid word  & $v_2$ & $v_3$ & $v_{4a}$ & $v_{4b}$ & $v_{5a}$ & $v_{5b}$ & $v_{5c}$ & $v_{6a}$ & $v_{6b}$ & $v_{6c}$ & $v_{6d}$ & $v_{6e}$ \\
\hline
\hline\newstrut
03.001 & aaa &  1 &  -1 &  1 &  -3 &  -3 &  1 &  -2 &  3 &  -3 &  2 &  -3 &  -1 \\ 
\hline\newstrut
04.001 & aBaB & -1 &  0 &  -2 &  3 &  0 &  0 &  0 &  -6 &  4 &  -4 &  2 &  1 \\ 
\hline\newstrut
05.001 & aaaaa &  3 &  -5 &  1 &  -6 & -12 &  4 &  -8 &  -6 &  -6 &  0 &  -7 &  2 \\ 
\hline\newstrut
05.002 & aabbAb &  2 &  -3 &  1 &  -5 &  -7 &  3 &  -5 &  1 &  -5 &  1 &  -6 &  0 \\ 
\hline\newstrut
06.001 & AbAcBcb & -2 &  1 &  -5 &  5 &  4 &  -1 &  2 & -10 &  7 &  -8 &  2 &  0 \\ 
\hline\newstrut
06.002 & AbAbbb & -1 &  1 &  -3 &  1 &  3 &  -1 &  1 &  -1 &  3 &  -3 &  0 &  -3 \\ 
\hline\newstrut
06.003 & AbbAAb &  1 &  0 &  2 &  -2 &  0 &  0 &  0 &  0 &  -4 &  2 &  -1 &  1 \\ 
\hline\newstrut
07.001 & aaaaaaa &  6 & -14 &  -4 &  -3 & -21 &  7 & -14 & -33 &  -7 &  -8 &  -5 &  8 \\ 
\hline\newstrut
07.002 & AcccbaaCb &  3 &  -6 &  0 &  -5 &  -9 &  6 &  -7 &  -3 &  -8 &  -3 &  -9 &  4 \\ 
\hline\newstrut
07.003 & aabAbbbb &  5 & -11 &  -3 &  -6 & -16 &  8 & -13 & -18 &  -8 &  -9 &  -9 &  7 \\ 
\hline\newstrut
07.004 & aabccAbCb &  4 &  -8 &  -2 &  -8 & -10 &  8 & -10 &  0 &  -8 &  -6 & -13 &  3 \\ 
\hline\newstrut
07.005 & aaaabAbb &  4 &  -8 &  0 &  -5 & -14 &  6 &  -9 & -12 &  -8 &  -3 &  -9 &  5 \\ 
\hline\newstrut
07.006 & aBAAcbbbc &  1 &  -2 &  0 &  -3 &  -2 &  3 &  -2 &  4 &  -3 &  0 &  -5 &  -1 \\ 
\hline\newstrut
07.007 & aCbCbAbCb & -1 &  1 &  -1 &  4 &  0 &  -2 &  0 &  -8 &  4 &  -4 &  4 &  3 \\ 
\hline\newstrut
08.001 & AbcBAddcbD & -3 &  3 &  -9 &  5 &  12 &  -3 &  5 &  -6 &  5 &  -9 &  0 &  -4 \\ 
\hline\newstrut
08.002 & AbbbbbAb &  0 &  1 &  -3 &  -6 &  2 &  0 &  -3 &  8 &  -4 &  1 &  -2 &  -7 \\ 
\hline\newstrut
08.003 & AABaddcDBc & -4 &  0 & -14 &  8 &  0 &  0 &  0 & -24 &  8 & -18 &  -3 &  -3 \\ 
\hline\newstrut
08.004 & aaacBCCaB & -3 &  1 & -11 &  4 &  0 &  -2 &  -1 & -14 &  6 & -11 &  -2 &  -6 \\ 
\hline\newstrut
08.005 & aaaBaaaB & -1 &  3 &  -5 &  -5 &  5 &  -3 &  -2 &  3 &  -3 &  0 &  2 & -10 \\ 
\hline\newstrut
08.006 & AbACbbbcc & -2 &  3 &  -7 &  0 &  9 &  -3 &  2 &  1 &  1 &  -4 &  0 &  -8 \\ 
\hline\newstrut
08.007 & aaaaBBaB &  2 &  -2 &  4 &  -2 &  -7 &  1 &  -3 &  -9 &  -7 &  1 &  -1 &  4 \\ 
\hline\newstrut
08.008 & AbaaCbbCC &  2 &  -1 &  3 &  -4 &  -2 &  1 &  -1 &  -4 &  -7 &  1 &  -3 &  1 \\ 
\hline\newstrut
08.009 & AbAAAbbb & -2 &  0 &  -8 &  1 &  0 &  0 &  0 &  -6 &  4 &  -6 &  -2 &  -7 \\ 
\hline\newstrut
08.010 & AbbAAbbb &  3 &  -3 &  3 &  -6 &  -5 &  3 &  -3 &  -3 &  -9 &  1 &  -7 &  1 \\ 
\hline\newstrut
08.011 & AbbCbccAb & -1 &  2 &  -4 &  -2 &  8 &  -1 &  2 &  10 &  -1 &  0 &  -2 &  -7 \\ 
\hline\newstrut
08.012 & aBcDcDbaCB & -3 &  0 &  -8 &  8 &  0 &  0 &  0 & -16 &  10 & -12 &  2 &  2 \\ 
\hline\newstrut
08.013 & aabCbACCb &  1 &  -1 &  3 &  0 &  -6 &  0 &  -3 & -10 &  -4 &  -1 &  2 &  4 \\ 
\hline\newstrut
08.014 & aabbACbCb &  0 &  0 &  -2 &  -3 &  -2 &  0 &  -3 &  4 &  0 &  1 &  -2 &  -4 \\ 
\hline\newstrut
08.015 & aaBaccbbc &  4 &  -7 &  1 &  -7 & -16 &  5 & -10 &  -8 &  -7 &  2 &  -9 &  3 \\ 
\hline\newstrut
08.016 & aaBaaBaB &  1 &  -1 &  3 &  0 &  2 &  2 &  2 &  -4 &  -4 &  -2 &  -1 &  2 \\ 
\hline\newstrut
08.017 & AbAbbAAb & -1 &  0 &  -4 &  0 &  0 &  0 &  0 &  2 &  2 &  0 &  -1 &  -3 \\ 
\hline\newstrut
08.018 & aBaBaBaB &  1 &  0 &  0 &  -5 &  0 &  0 &  0 &  14 &  -4 &  8 &  -5 &  -2 \\ 
\hline\newstrut
08.019 & abababba &  5 & -10 &  0 &  -5 & -18 &  6 & -10 & -14 &  -8 &  0 & -11 &  4 \\ 
\hline\newstrut
08.020 & aaabAAAb &  2 &  -2 &  2 &  -5 &  -1 &  3 &  -1 &  5 &  -7 &  3 &  -6 &  0 \\ 
\hline\newstrut
08.021 & aBBaabbb &  0 &  1 &  -1 &  -3 &  1 &  -1 &  -1 &  5 &  -1 &  1 &  -1 &  -5 \\ 
\hline\newstrut
09.001 & aaaaaaaaa & 10 & -30 & -20 &  15 &  -9 &  3 &  -6 & -27 & -13 &  0 &  21 & -12 \\ 
\hline\newstrut
09.002 & abcdddcDbaCB &  4 & -10 &  -2 &  -2 &  -5 &  10 &  -5 &  1 & -16 &  -9 & -14 &  10 \\ 
\hline\newstrut
09.003 & aBaaaaaabb &  9 & -26 & -18 &  6 &  -7 &  9 & -12 & -23 & -15 & -18 &  11 &  5 \\ 
\hline\newstrut
09.004 & AcbaabccccB &  7 & -19 & -11 &  -1 &  -5 &  13 & -11 &  -5 & -19 & -21 &  -9 &  12 \\ 
\hline\newstrut
09.005 & aabAcBcddcDb &  6 & -15 &  -9 &  -7 &  0 &  15 & -10 &  20 & -17 & -18 & -20 &  10 \\ 
\hline\newstrut
09.006 & aabbaaaaaB &  7 & -18 &  -6 &  3 & -17 &  7 &  -8 & -25 & -11 &  -4 &  -8 &  3 \\ 
\hline\newstrut
09.007 & aaabccAbbbC &  5 & -12 &  -2 &  -1 &  -9 &  9 &  -5 &  -5 & -15 &  -7 & -15 &  7 \\ 
\hline\newstrut
09.008 & aBcAADcBccdc &  0 &  -2 &  -2 &  -1 &  1 &  4 &  0 &  3 &  0 &  -2 &  -4 &  -2 \\ 
\hline\newstrut
09.009 & aaaBaaaabb &  8 & -22 & -12 &  4 & -11 &  9 & -10 & -23 & -15 & -14 &  -1 &  5 \\ 
\hline\newstrut
09.010 & AbaabbbCbcc &  8 & -22 & -16 &  -2 &  -2 &  14 & -14 &  0 & -16 & -26 &  -4 &  12 \\ 
\hline\newstrut
09.011 & AbCbAbbbbcb &  4 &  -9 &  -3 &  -7 &  -6 &  10 &  -8 &  4 & -10 &  -8 & -15 &  4 \\ 
\hline\newstrut
09.012 & aBAAcbbbddcD &  1 &  -3 &  -1 &  -2 &  4 &  6 &  1 &  10 &  -6 &  -3 &  -9 &  1 \\ 
\hline\newstrut
09.013 & aabCbAccbbb &  7 & -18 & -10 &  -3 &  -8 &  12 & -12 &  -8 & -14 & -18 & -10 &  9 \\ 
\hline\newstrut
09.014 & addCbCbCADbbcB & -1 &  2 &  0 &  4 &  -4 &  -5 &  -2 & -12 &  5 &  -4 &  7 &  2 \\ 
\hline\newstrut
09.015 & aBacBdCddc &  2 &  -5 &  -1 &  -4 &  -3 &  7 &  -4 &  5 &  -7 &  -4 & -10 &  2 \\ 
\hline\newstrut
09.016 & aaBaaabbbb &  6 & -14 &  -2 &  0 & -12 &  8 &  -4 &  -2 & -14 &  -2 & -19 &  1 \\ 
\hline\newstrut
09.017 & aaacBaBCBaB & -2 &  0 &  -4 &  7 &  5 &  1 &  4 &  -9 &  9 & -10 &  1 &  2 \\ 
\hline\newstrut
09.018 & aaCbAbbccbb &  6 & -15 &  -5 &  -1 & -14 &  9 & -10 & -20 & -13 & -10 &  -9 &  9 \\ 
\hline
\end{tabular}}
\newpage

\cl{\begin{tabular}{|l|l|r|r|rr|rrr|rrrrr|} \hline\newstrut
Knot & Braid word  & $v_2$ & $v_3$ & $v_{4a}$ & $v_{4b}$ & $v_{5a}$ & $v_{5b}$ & $v_{5c}$ & $v_{6a}$ & $v_{6b}$ & $v_{6c}$ & $v_{6d}$ & $v_{6e}$ \\
\hline
\hline\newstrut
09.019 & aabACDCbCbdC & -2 &  -1 &  -3 &  8 &  4 &  3 &  3 & -16 &  9 & -11 &  5 &  5 \\ 
\hline\newstrut
09.020 & abbCbAbCbbb &  2 &  -4 &  0 &  -4 &  0 &  6 &  -1 &  8 &  -8 &  -3 & -11 &  0 \\ 
\hline\newstrut
09.021 & AAcDcBacddbb &  3 &  -6 &  -2 &  -8 &  0 &  9 &  -4 &  20 &  -9 &  -2 & -17 &  0 \\ 
\hline\newstrut
09.022 & aBcccBcABcB & -1 &  -1 &  -1 &  5 &  5 &  3 &  4 &  -7 &  5 &  -8 &  1 &  2 \\ 
\hline\newstrut
09.023 & aBaabcccbCb &  5 & -11 &  -1 &  -3 & -15 &  7 &  -8 & -15 &  -9 &  -4 & -12 &  4 \\ 
\hline\newstrut
09.024 & accBacBBB &  1 &  -2 &  0 &  -4 &  2 &  4 &  1 &  8 &  -4 &  3 &  -5 &  -2 \\ 
\hline\newstrut
09.025 & AbADCbddcddbbC &  0 &  -1 &  -3 &  -2 &  0 &  2 &  -2 &  0 &  0 &  -4 &  -4 &  -3 \\ 
\hline\newstrut
09.026 & AbbbcbAbACb &  0 &  1 &  1 &  0 &  -5 &  -3 &  -3 &  -3 &  1 &  1 &  3 &  0 \\ 
\hline\newstrut
09.027 & AbAACbAbbcb &  0 &  -1 &  -1 &  -1 &  1 &  2 &  1 &  -1 &  0 &  -1 &  -1 &  -2 \\ 
\hline\newstrut
09.028 & aaccBBacB &  1 &  0 &  2 &  -3 &  -1 &  0 &  -1 &  5 &  -4 &  3 &  -2 &  -1 \\ 
\hline\newstrut
09.029 & AbCbAbCbb &  1 &  -2 &  2 &  -1 &  10 &  6 &  6 &  0 &  -4 &  -4 &  -4 &  -1 \\ 
\hline\newstrut
09.030 & aCCbCAbaaCb & -1 &  -1 &  -3 &  2 &  -3 &  1 &  -1 &  -5 &  3 &  -3 &  1 &  0 \\ 
\hline\newstrut
09.031 & AbbCbbAbcAb &  2 &  -2 &  2 &  -5 &  -9 &  1 &  -6 &  1 &  -5 &  4 &  -4 &  1 \\ 
\hline\newstrut
09.032 & AbcAbAccbCb & -1 &  2 &  -2 &  1 &  2 &  -3 &  0 &  0 &  3 &  -2 &  2 &  -3 \\ 
\hline\newstrut
09.033 & AAAbAbbcaBc &  1 &  -1 &  1 &  -3 &  5 &  3 &  3 &  7 &  -5 &  5 &  -3 &  1 \\ 
\hline\newstrut
09.034 & aBcBaBcaB & -1 &  0 &  -2 &  3 &  8 &  2 &  5 &  -2 &  2 &  -1 &  2 &  3 \\ 
\hline\newstrut
09.035 & accDBabccbCdCb &  7 & -18 & -14 &  -9 &  6 &  18 & -12 &  38 & -16 & -24 & -22 &  11 \\ 
\hline\newstrut
09.036 & AbbbCbccAbc &  3 &  -7 &  -1 &  -4 &  -9 &  7 &  -7 &  -9 &  -7 &  -7 &  -9 &  4 \\ 
\hline\newstrut
09.037 & aBcBcAdCBdcB & -3 &  -1 &  -7 &  10 &  0 &  2 &  1 & -20 &  12 & -15 &  3 &  5 \\ 
\hline\newstrut
09.038 & aabccbbAbCb &  6 & -14 &  -4 &  -3 & -13 &  9 &  -9 &  -9 & -13 &  -5 & -13 &  7 \\ 
\hline\newstrut
09.039 & ACbdccAbbcdB &  2 &  -4 &  -2 &  -7 &  1 &  7 &  -3 &  19 &  -7 &  -1 & -14 &  -2 \\ 
\hline\newstrut
09.040 & acBcaBacB & -1 &  1 &  -3 &  1 &  -5 &  -3 &  -4 &  -5 &  5 &  -2 &  2 &  -2 \\ 
\hline\newstrut
09.041 & AbbDcDBaCbcDcb &  0 &  1 &  3 &  3 & -11 &  -5 &  -5 & -23 &  3 &  -5 &  9 &  3 \\ 
\hline\newstrut
09.042 & aaacBcAAB & -2 &  0 &  -6 &  4 &  0 &  0 &  0 &  -6 &  6 &  -6 &  0 &  -1 \\ 
\hline\newstrut
09.043 & abaabbcBaBC &  1 &  -2 &  -2 &  -5 &  -7 &  2 &  -7 &  -3 &  -2 &  -1 &  -3 &  -1 \\ 
\hline\newstrut
09.044 & AbAcBcbbC &  0 &  -1 &  1 &  2 &  1 &  2 &  1 &  -7 &  0 &  -3 &  2 &  3 \\ 
\hline\newstrut
09.045 & aBacbbbcB &  2 &  -4 &  0 &  -4 &  -8 &  4 &  -6 &  -6 &  -4 &  -4 &  -6 &  1 \\ 
\hline\newstrut
09.046 & acbACbacB & -2 &  3 &  -5 &  3 &  11 &  -3 &  5 &  3 &  3 &  -3 &  1 &  -3 \\ 
\hline\newstrut
09.047 & AbcAbAbcb & -1 &  2 &  0 &  4 &  4 &  -3 &  3 &  -6 &  5 &  -5 &  4 &  0 \\ 
\hline\newstrut
09.048 & aaBcbbACbCb &  3 &  -5 &  -1 &  -9 &  -3 &  7 &  -5 &  17 &  -7 &  1 & -15 &  -2 \\ 
\hline\newstrut
09.049 & aabbcbAbbcB &  6 & -14 &  -6 &  -6 & -14 &  10 & -14 & -14 & -10 & -12 & -10 &  9 \\ 
\hline\newstrut
10.001 & aaDECbAcdEdcb & -4 &  6 & -14 &  2 &  23 &  -6 &  7 &  11 &  -8 &  -3 &  -2 & -10 \\ 
\hline\newstrut
10.002 & aBaaaaaaaB &  2 &  -2 &  -4 & -15 &  -7 &  5 & -16 &  -9 & -15 &  -4 &  4 &  5 \\ 
\hline\newstrut
10.003 & aCBaEEbcDecbD & -6 &  3 & -27 &  7 &  15 &  -3 &  5 & -29 & -13 & -21 & -16 & -10 \\ 
\hline\newstrut
10.004 & aaDDDbcBADcb & -5 &  1 & -23 &  4 &  16 &  2 &  8 & -22 &  -8 & -14 & -12 & -11 \\ 
\hline\newstrut
10.005 & AAbAbbbbbb &  4 &  -7 &  5 &  0 & -20 &  2 &  -7 & -32 &  -8 &  -3 &  -4 &  5 \\ 
\hline\newstrut
10.006 & aBcABAcbbbbbb & -1 &  4 &  -8 & -12 &  7 &  -2 &  -8 &  5 & -16 &  2 &  7 & -10 \\ 
\hline\newstrut
10.007 & aCBabcDcDccb & -1 &  3 &  -5 &  -6 &  14 &  0 &  2 &  28 & -10 &  6 &  -3 &  -8 \\ 
\hline\newstrut
10.008 & aaaaacBaCCB & -3 &  4 & -16 &  -6 &  6 &  -5 &  -5 &  -8 & -11 &  -1 &  0 & -16 \\ 
\hline\newstrut
10.009 & AAAbAbbbbb & -2 &  2 & -12 &  -6 &  7 &  -1 &  -1 &  -1 &  -7 &  1 &  -2 & -13 \\ 
\hline\newstrut
10.010 & ACdccbaaCbDCbC &  1 &  -2 &  4 &  3 & -14 &  -1 &  -6 & -30 &  -1 &  -8 &  6 &  6 \\ 
\hline\newstrut
10.011 & AACBadcccBcd & -5 &  4 & -24 &  1 &  9 &  -6 &  -1 & -25 & -14 & -13 &  -9 & -16 \\ 
\hline\newstrut
10.012 & AAbAccbbbbC &  4 &  -6 &  4 &  -3 & -15 &  3 &  -7 & -33 & -11 &  -7 &  -3 &  5 \\ 
\hline\newstrut
10.013 & aaCbACdeedcEBBD & -5 &  2 & -18 &  10 &  6 &  -3 &  2 & -22 &  5 & -18 &  -6 &  -1 \\ 
\hline\newstrut
10.014 & AbAccbbbbCb &  2 &  -3 &  -3 & -11 & -10 &  4 & -13 &  0 &  -6 &  1 &  -3 &  1 \\ 
\hline\newstrut
10.015 & AcBccccBAAb &  3 &  -2 &  6 &  -2 & -10 &  0 &  -5 & -28 & -12 &  -5 &  1 &  6 \\ 
\hline\newstrut
10.016 & aBAAcBDcddcc & -4 &  4 & -18 &  0 &  15 &  -4 &  3 &  -7 &  -8 &  -7 &  -6 & -14 \\ 
\hline\newstrut
10.017 & aaaaBaBBBB &  2 &  0 &  8 &  3 &  0 &  0 &  0 & -24 & -10 &  -6 &  4 &  9 \\ 
\hline\newstrut
10.018 & accddCBdAAcB & -2 &  1 &  -9 &  -1 &  -6 &  -3 &  -6 &  -2 &  3 &  0 &  -1 &  -7 \\ 
\hline\newstrut
10.019 & AABaccccBcB &  1 &  0 &  6 &  4 &  5 &  1 &  4 & -21 &  -5 &  -8 &  5 &  7 \\ 
\hline\newstrut
10.020 & ACbbbcddcDAb & -3 &  6 & -12 &  -4 &  17 &  -6 &  1 &  9 & -12 &  1 &  3 & -13 \\ 
\hline\newstrut
10.021 & AbAccbCbbbb &  1 &  0 &  -4 & -13 &  3 &  4 &  -8 &  19 & -14 &  2 &  -4 &  -4 \\ 
\hline\newstrut
10.022 & AAAbccAbbbC & -4 &  2 & -20 &  -1 &  11 &  -1 &  3 & -17 &  -9 &  -9 &  -8 & -16 \\ 
\hline\newstrut
10.023 & AbAbbbCbAcc &  3 &  -5 &  5 &  0 & -19 &  1 &  -8 & -35 &  -7 &  -6 &  1 &  6 \\ 
\hline\newstrut
10.024 & AbCbbaabcDcD & -2 &  5 &  -9 &  -7 &  17 &  -3 &  1 &  21 & -13 &  3 &  0 & -14 \\ 
\hline\newstrut
10.025 & AccbbbCbbAb &  0 &  2 &  -4 & -10 &  7 &  1 &  -4 &  17 & -13 &  4 &  0 &  -6 \\ 
\hline
\end{tabular}}
\newpage

\cl{\begin{tabular}{|l|l|r|r|rr|rrr|rrrrr|} \hline\newstrut
Knot & Braid word  & $v_2$ & $v_3$ & $v_{4a}$ & $v_{4b}$ & $v_{5a}$ & $v_{5b}$ & $v_{5c}$ & $v_{6a}$ & $v_{6b}$ & $v_{6c}$ & $v_{6d}$ & $v_{6e}$ \\
\hline
\hline\newstrut
10.026 & abbCBBBaBaacB & -3 &  2 & -14 &  -1 &  14 &  0 &  5 &  -2 &  -4 &  -3 &  -5 & -11 \\ 
\hline\newstrut
10.027 & AAbCbccccAb &  2 &  -3 &  5 &  1 & -14 &  0 &  -5 & -22 &  -4 &  -3 &  1 &  5 \\ 
\hline\newstrut
10.028 & aabCbADcDDcb &  3 &  -4 &  4 &  -3 & -13 &  2 &  -7 & -27 & -10 &  -7 &  -1 &  4 \\ 
\hline\newstrut
10.029 & aBcDcBABcBdccc & -4 &  3 & -15 &  5 &  5 &  -5 &  0 & -13 &  3 & -10 &  -3 &  -6 \\ 
\hline\newstrut
10.030 & aabCbAccdCbCDb &  1 &  -1 &  -3 &  -9 &  -1 &  3 &  -6 &  17 &  -5 &  4 &  -7 &  -5 \\ 
\hline\newstrut
10.031 & aabbAdCDDCbC &  2 &  1 &  5 &  -1 &  10 &  1 &  6 & -22 &  -9 &  -6 &  2 &  4 \\ 
\hline\newstrut
10.032 & AbCbaaCCCbb & -1 &  0 &  -6 &  -3 &  -5 &  -1 &  -4 &  3 &  1 &  2 &  -2 &  -7 \\ 
\hline\newstrut
10.033 & aBcDcBAABcdd &  0 &  0 &  4 &  6 &  0 &  0 &  0 & -24 &  0 & -10 &  8 &  8 \\ 
\hline\newstrut
10.034 & ADDccbaaDCbc &  3 &  -3 &  3 &  -5 &  -4 &  3 &  -2 & -12 &  -9 &  -4 &  -6 &  1 \\ 
\hline\newstrut
10.035 & AbCbCdEdEaBcDbb & -4 &  2 & -12 &  9 &  9 &  -2 &  4 & -11 &  8 & -12 &  -1 &  1 \\ 
\hline\newstrut
10.036 & aaCbbcDcDAcb &  1 &  -2 &  -2 &  -6 &  -9 &  2 &  -9 &  1 &  -2 &  3 &  -2 &  0 \\ 
\hline\newstrut
10.037 & adCDDbCCbaaB &  3 &  0 &  4 &  -5 &  0 &  0 &  0 & -18 & -10 &  -4 &  -3 &  2 \\ 
\hline\newstrut
10.038 & aadCBBcccddCAb & -1 &  2 &  -6 &  -5 &  1 &  -2 &  -5 &  1 &  -2 &  1 &  2 &  -7 \\ 
\hline\newstrut
10.039 & AbAbbbCbccc &  1 &  -1 &  -3 &  -9 &  -6 &  2 & -10 &  -4 &  -6 &  0 &  1 &  -1 \\ 
\hline\newstrut
10.040 & aabAbbCbaCC &  3 &  -4 &  4 &  -3 &  -8 &  3 &  -3 & -10 &  -9 &  -1 &  -6 &  3 \\ 
\hline\newstrut
10.041 & aBcBADcBcdcBcc & -2 &  2 &  -6 &  2 &  10 &  -1 &  4 &  6 &  3 &  -2 &  -1 &  -4 \\ 
\hline\newstrut
10.042 & aBAAcDbbcDDc &  0 &  -1 &  1 &  2 &  -7 &  0 &  -4 & -13 &  0 &  -6 &  3 &  2 \\ 
\hline\newstrut
10.043 & AbCDbCCabCDcbb &  2 &  0 &  2 &  -5 &  0 &  0 &  0 &  -4 &  -6 &  -2 &  -5 &  -2 \\ 
\hline\newstrut
10.044 & ABcBcddabCBBdc &  0 &  -1 &  -1 &  -1 &  -4 &  1 &  -3 &  4 &  1 &  3 &  -1 &  0 \\ 
\hline\newstrut
10.045 & aCbCbdCAbCbDCb & -2 &  0 &  -4 &  7 &  0 &  0 &  0 & -18 &  8 & -12 &  4 &  3 \\ 
\hline\newstrut
10.046 & AbbbbbAbbb &  0 &  4 &  -6 & -18 &  -2 &  -2 & -17 & -18 & -20 &  1 &  20 &  -5 \\ 
\hline\newstrut
10.047 & aaaaaBaaBB &  6 & -11 &  -1 &  -6 & -10 &  8 &  -7 & -16 & -12 &  -9 & -14 &  3 \\ 
\hline\newstrut
10.048 & AAbbbbAAAb &  4 &  0 &  6 &  -5 & -11 &  -3 &  -6 & -21 & -15 &  -4 &  -5 &  4 \\ 
\hline\newstrut
10.049 & aaaacbbbcaB &  7 & -16 &  -4 &  -3 & -23 &  7 & -13 & -21 &  -7 &  3 & -10 &  4 \\ 
\hline\newstrut
10.050 & AbbCbccAbbb & -1 &  5 &  -7 & -13 &  9 &  -3 &  -7 &  7 & -17 &  3 &  9 & -13 \\ 
\hline\newstrut
10.051 & aabCCAbbCbb &  5 &  -8 &  0 &  -8 &  -5 &  8 &  -5 &  -5 & -12 &  -7 & -15 &  2 \\ 
\hline\newstrut
10.052 & AABccBcccaB &  3 &  -1 &  5 &  -4 &  10 &  4 &  6 & -14 & -12 &  -6 &  -5 &  1 \\ 
\hline\newstrut
10.053 & aBccdbAccBcAddbb &  6 & -13 &  -3 &  -6 & -19 &  8 & -13 & -15 &  -8 &  -1 & -11 &  6 \\ 
\hline\newstrut
10.054 & AABccBcccAb &  4 &  -2 &  4 &  -7 &  -2 &  2 &  -1 & -16 & -12 &  -3 &  -7 &  2 \\ 
\hline\newstrut
10.055 & AdcdbaaabccDcB &  5 & -10 &  0 &  -6 & -17 &  7 & -10 &  -7 &  -9 &  2 & -13 &  4 \\ 
\hline\newstrut
10.056 & AbbccBcAbbb &  0 &  2 &  -4 & -10 &  -4 &  -2 & -10 &  -4 &  -6 &  4 &  7 &  -6 \\ 
\hline\newstrut
10.057 & aaabCCAbbCb &  4 &  -6 &  2 &  -6 &  -6 &  6 &  -4 &  -2 & -12 &  -2 & -12 &  3 \\ 
\hline\newstrut
10.058 & abEdCbCDEaBcd & -4 &  1 & -11 &  11 &  5 &  -1 &  3 & -21 &  13 & -17 &  1 &  3 \\ 
\hline\newstrut
10.059 & abbdAACdCBcddC & -1 &  1 &  -3 &  0 &  2 &  -1 &  0 &  4 &  3 &  0 &  0 &  -4 \\ 
\hline\newstrut
10.060 & AAbbaBcDcBcBdc & -1 &  1 &  -1 &  3 &  4 &  -1 &  3 &  -4 &  3 &  1 &  5 &  3 \\ 
\hline\newstrut
10.061 & aBAAcccBccc & -4 &  5 & -21 &  -5 &  1 &  -9 &  -9 & -29 & -15 &  -5 &  2 & -18 \\ 
\hline\newstrut
10.062 & AAbbbAbbbb &  5 &  -9 &  3 &  -2 & -15 &  5 &  -6 & -29 & -13 &  -6 &  -8 &  5 \\ 
\hline\newstrut
10.063 & aabccAdCBcddcb &  6 & -14 &  -4 &  -4 & -20 &  8 & -14 & -24 &  -8 &  -4 &  -7 &  9 \\ 
\hline\newstrut
10.064 & AAAbbbAbbb & -3 &  3 & -17 &  -6 &  4 &  -4 &  -4 & -18 & -10 &  -4 &  -1 & -18 \\ 
\hline\newstrut
10.065 & acbbcABcccABB &  4 &  -7 &  3 &  -3 & -13 &  5 &  -7 & -23 & -13 &  -7 &  -6 &  6 \\ 
\hline\newstrut
10.066 & aabCbbaBccccb &  7 & -17 &  -5 &  0 & -24 &  6 & -13 & -34 &  -6 &  1 &  -4 &  6 \\ 
\hline\newstrut
10.067 & AbDCbDaaccbbdC &  0 &  0 &  -4 &  -6 &  2 &  2 &  -4 &  8 &  -4 &  0 &  -3 &  -5 \\ 
\hline\newstrut
10.068 & AbCddaaBcbbDCC &  2 &  -3 &  5 &  1 & -19 &  -1 &  -9 & -39 &  -5 &  -9 &  6 &  6 \\ 
\hline\newstrut
10.069 & acdbCdBBAAcbbb &  2 &  -4 &  0 &  -4 &  -5 &  5 &  -5 &  5 &  -7 &  -3 &  -9 &  2 \\ 
\hline\newstrut
10.070 & aBcccAAbbDcD & -3 &  2 & -10 &  4 &  7 &  -2 &  2 &  -1 &  4 &  -6 &  -2 &  -4 \\ 
\hline\newstrut
10.071 & aaDDBCbbCDbcAb &  1 &  0 &  2 &  -1 &  0 &  0 &  0 &  -8 &  -4 &  -2 &  0 &  2 \\ 
\hline\newstrut
10.072 & aaaBaabCbCb &  2 &  -4 &  -2 &  -7 & -15 &  3 & -13 & -13 &  -3 &  1 &  0 &  5 \\ 
\hline\newstrut
10.073 & aBccaDcBcdcB &  1 &  -2 &  0 &  -2 &  -5 &  2 &  -4 &  -3 &  -2 &  -4 &  -4 &  0 \\ 
\hline\newstrut
10.074 & AdbCbDCbaaDccb &  0 &  2 &  -4 & -10 &  12 &  2 &  0 &  38 & -12 &  8 &  -8 & -10 \\ 
\hline\newstrut
10.075 & aBacBcdcBcDB &  0 &  1 &  1 &  0 &  0 &  -2 &  1 &  -4 &  0 &  3 &  4 &  1 \\ 
\hline\newstrut
10.076 & aaBAcccBaBccc & -2 &  6 & -10 & -10 &  10 &  -6 &  -6 &  0 & -16 &  4 &  11 & -13 \\ 
\hline\newstrut
10.077 & AbbbCCaabbC &  4 &  -5 &  3 &  -6 &  -8 &  4 &  -4 & -14 & -10 &  -4 &  -9 &  1 \\ 
\hline\newstrut
10.078 & AAdcBaccDcbbbc &  3 &  -5 &  1 &  -7 & -10 &  5 &  -7 &  2 &  -7 &  3 &  -9 &  1 \\ 
\hline\newstrut
10.079 & aaaBBaaBBB &  5 &  0 &  4 &  -9 &  0 &  0 &  0 & -16 & -14 &  -2 & -10 &  3 \\ 
\hline\newstrut
10.080 & aabCbbAbbbcbb &  6 & -12 &  0 &  -6 & -21 &  7 & -11 &  -7 &  -9 &  5 & -16 &  2 \\ 
\hline\newstrut
10.081 & AABBACbddccd &  3 &  0 &  2 &  -8 &  0 &  0 &  0 &  8 & -10 &  4 & -10 &  0 \\ 
\hline
\end{tabular}}
\newpage

\cl{\begin{tabular}{|l|l|r|r|rr|rrr|rrrrr|} \hline\newstrut
Knot & Braid word  & $v_2$ & $v_3$ & $v_{4a}$ & $v_{4b}$ & $v_{5a}$ & $v_{5b}$ & $v_{5c}$ & $v_{6a}$ & $v_{6b}$ & $v_{6c}$ & $v_{6d}$ & $v_{6e}$ \\
\hline
\hline\newstrut
10.082 & aaaaBBaBaB &  0 &  0 &  -4 &  -6 &  1 &  1 &  -2 &  15 &  -5 &  10 &  -3 &  -3 \\ 
\hline\newstrut
10.083 & aabCbbCbACb &  1 &  -2 &  4 &  3 &  2 &  3 &  4 &  -8 &  -3 &  -8 &  -3 &  2 \\ 
\hline\newstrut
10.084 & aaabbCbACbC &  2 &  -2 &  4 &  -2 &  1 &  3 &  2 &  -3 &  -7 &  -4 &  -6 &  -1 \\ 
\hline\newstrut
10.085 & AbbAbAbbbb &  2 &  -3 &  5 &  1 &  2 &  4 &  5 &  0 &  -6 &  -7 & -10 &  -2 \\ 
\hline\newstrut
10.086 & AbCbCCbaabC & -1 &  1 &  -5 &  -2 &  6 &  0 &  2 &  14 &  -2 &  8 &  -2 &  -2 \\ 
\hline\newstrut
10.087 & aaaCbCAbCbC &  0 &  -1 &  -3 &  -4 &  -4 &  1 &  -3 &  8 &  -1 &  5 &  -3 &  -4 \\ 
\hline\newstrut
10.088 & aaBDcaBcBdCAbC & -1 &  0 &  0 &  6 &  0 &  0 &  0 & -16 &  4 &  -8 &  6 &  6 \\ 
\hline\newstrut
10.089 & abCbAdbCbdCd &  1 &  -3 &  -1 &  -2 &  -4 &  4 &  -4 &  -4 &  -2 &  -8 &  -6 &  -1 \\ 
\hline\newstrut
10.090 & aacBBaaBaBcAB & -3 &  1 & -13 &  1 &  3 &  -1 &  0 &  -5 &  -1 &  -2 &  -3 &  -6 \\ 
\hline\newstrut
10.091 & aaaBBaBBaB &  2 &  0 &  6 &  0 & -11 &  -3 &  -6 & -19 &  -7 &  -8 &  0 &  2 \\ 
\hline\newstrut
10.092 & aaabbCbAbCb &  2 &  -3 &  -1 &  -8 & -14 &  2 & -12 &  -6 &  -4 &  6 &  0 &  4 \\ 
\hline\newstrut
10.093 & AAccBcaBccB &  1 &  -1 &  5 &  3 &  13 &  5 &  8 & -25 &  -3 & -14 &  3 &  2 \\ 
\hline\newstrut
10.094 & aaaBaaBBaB & -2 &  2 & -10 &  -3 &  1 &  -3 &  -3 &  -3 &  -3 &  3 &  1 &  -7 \\ 
\hline\newstrut
10.095 & AAbCbAbbccb &  3 &  -5 &  3 &  -3 & -10 &  4 &  -5 & -12 &  -8 &  -5 &  -7 &  2 \\ 
\hline\newstrut
10.096 & abCdbacBCCdC & -3 &  2 &  -8 &  7 &  9 &  -2 &  5 &  -3 &  6 &  -3 &  2 &  4 \\ 
\hline\newstrut
10.097 & aaBcDbbACbccDb &  2 &  -4 &  -4 & -10 &  -6 &  6 & -10 &  14 &  -4 &  0 & -11 &  -2 \\ 
\hline\newstrut
10.098 & AbbccbAbbCb &  0 &  3 &  -3 & -11 &  5 &  -1 &  -5 &  19 & -13 &  9 &  3 &  -7 \\ 
\hline\newstrut
10.099 & AAbAAbbAbb &  4 &  0 &  4 &  -8 &  0 &  0 &  0 & -10 & -10 &  -4 & -10 &  -2 \\ 
\hline\newstrut
10.100 & aaaBaaBaaB &  4 &  -7 &  3 &  -3 &  3 &  9 &  3 &  -1 & -15 & -11 & -17 &  -1 \\ 
\hline\newstrut
10.101 & AcbbadcbbbccdCAB &  7 & -17 &  -7 &  -3 & -17 &  9 & -13 & -19 &  -9 &  -5 &  -8 &  7 \\ 
\hline\newstrut
10.102 & abbCbACbaCC & -2 &  1 &  -9 &  -1 &  10 &  1 &  4 &  6 &  -1 &  2 &  -3 &  -6 \\ 
\hline\newstrut
10.103 & aabbCbbCAbC &  3 &  -4 &  4 &  -3 &  3 &  6 &  3 &  -5 & -12 &  -7 & -10 &  0 \\ 
\hline\newstrut
10.104 & aaBBBaaBaB &  1 &  0 &  6 &  4 &  11 &  3 &  6 & -19 &  -3 &  -8 &  4 &  6 \\ 
\hline\newstrut
10.105 & ABcBdacBccddCB & -1 &  0 &  -4 &  0 &  -5 &  -1 &  -4 &  7 &  3 &  4 &  -1 &  -1 \\ 
\hline\newstrut
10.106 & aaaBBaaBaB & -1 &  1 &  -7 &  -5 &  4 &  0 &  -1 &  8 &  -4 &  5 &  -2 &  -7 \\ 
\hline\newstrut
10.107 & adcBccdCBBABcB &  1 &  -1 &  3 &  0 &  -3 &  1 &  -2 & -13 &  -5 &  -4 &  2 &  3 \\ 
\hline\newstrut
10.108 & AAccBaccBcB &  0 &  0 &  4 &  6 &  16 &  4 &  10 & -22 &  2 & -14 &  5 &  4 \\ 
\hline\newstrut
10.109 & AAbbAAbbAb &  3 &  0 &  6 &  -3 &  0 &  0 &  0 & -14 & -12 &  -2 &  -3 &  4 \\ 
\hline\newstrut
10.110 & aCbCDcbbbaBCDb & -3 &  3 & -11 &  2 &  7 &  -4 &  1 &  1 &  2 &  -3 &  -2 &  -7 \\ 
\hline\newstrut
10.111 & aabbCbbAbCb &  1 &  0 &  -2 & -10 & -12 &  -1 & -13 & -10 &  -3 &  5 &  6 &  -2 \\ 
\hline\newstrut
10.112 & aaaBaBaBaB &  2 &  -2 &  0 &  -8 &  -6 &  2 &  -5 &  18 &  -8 &  15 &  -6 &  3 \\ 
\hline\newstrut
10.113 & aaabCbAbCbC &  0 &  1 &  1 &  0 &  3 &  -1 &  2 &  3 &  1 &  -4 &  -2 &  -5 \\ 
\hline\newstrut
10.114 & ABcBcBBaaCbcc &  1 &  -1 &  -1 &  -6 &  0 &  2 &  -1 &  24 &  -6 &  15 &  -6 &  1 \\ 
\hline\newstrut
10.115 & AbcDcBAdcdBcBB &  1 &  0 &  4 &  1 &  0 &  0 &  0 & -16 &  -6 &  -4 &  4 &  5 \\ 
\hline\newstrut
10.116 & aaBaBaBaaB &  0 &  0 &  -2 &  -3 &  -5 &  -1 &  -4 &  7 &  -3 &  10 &  1 &  2 \\ 
\hline\newstrut
10.117 & AbCbaabCbbC &  2 &  -3 &  3 &  -2 &  -5 &  3 &  -2 &  -7 &  -5 &  -6 &  -6 &  -1 \\ 
\hline\newstrut
10.118 & aaBaBBaBaB &  0 &  0 &  2 &  3 &  0 &  0 &  0 & -14 &  4 & -12 &  1 &  -2 \\ 
\hline\newstrut
10.119 & aabCCbACbCb & -1 &  0 &  -4 &  0 &  3 &  1 &  1 &  9 &  -1 &  7 &  0 &  2 \\ 
\hline\newstrut
10.120 & abbcDacBcAddccbC &  6 & -13 &  -1 &  -3 & -18 &  7 &  -8 &  0 & -13 &  12 & -15 &  7 \\ 
\hline\newstrut
10.121 & AbCbCbaabCb &  1 &  -2 &  2 &  0 &  0 &  3 &  1 &  -6 &  -1 &  -9 &  -5 &  -3 \\ 
\hline\newstrut
10.122 & aBcBcBABcbAbb &  2 &  -2 &  0 &  -8 &  2 &  4 &  0 &  24 &  -8 &  14 &  -9 &  1 \\ 
\hline\newstrut
10.123 & aBaBaBaBaB & -2 &  0 &  -6 &  4 &  0 &  0 &  0 & -14 &  10 & -16 &  -2 &  -8 \\ 
\hline\newstrut
10.124 & abbbbbabbb &  8 & -20 &  -6 &  4 & -18 &  6 &  -5 &  -4 & -10 &  9 & -17 &  -8 \\ 
\hline\newstrut
10.125 & ABBBAbbbbb &  3 &  0 &  4 &  -6 & -11 &  -3 &  -6 &  -3 & -11 &  4 &  -5 &  3 \\ 
\hline\newstrut
10.126 & aBBBabbbbb &  5 &  -9 &  -1 &  -8 &  -2 &  10 &  -4 &  6 & -14 &  -4 & -17 &  4 \\ 
\hline\newstrut
10.127 & aaaaabAAbb &  1 &  1 &  -1 & -11 &  -4 &  0 &  -9 &  0 &  -8 &  1 &  3 &  -6 \\ 
\hline\newstrut
10.128 & aacbbcbbAcb &  7 & -17 &  -5 &  0 & -16 &  8 &  -8 & -10 & -12 &  0 & -14 &  2 \\ 
\hline\newstrut
10.129 & AAcBBcbbacB &  2 &  -1 &  3 &  -4 &  9 &  4 &  5 &  -1 &  -8 &  1 &  -4 &  1 \\ 
\hline\newstrut
10.130 & AbCCCbaabcc &  4 &  -6 &  0 &  -9 &  3 &  9 &  -1 &  17 & -13 &  1 & -17 &  2 \\ 
\hline\newstrut
10.131 & AbcccbaabCC &  0 &  2 &  -2 &  -7 &  3 &  -1 &  -3 &  9 &  -5 &  1 &  0 &  -9 \\ 
\hline\newstrut
10.132 & AbcaaaBBBcb &  3 &  -5 &  1 &  -6 &  -1 &  7 &  -2 &  9 & -11 &  0 & -12 &  2 \\ 
\hline\newstrut
10.133 & AcbccbaaCBaCb &  1 &  0 &  0 &  -7 &  -5 &  0 &  -6 &  5 &  -4 &  4 &  -1 &  -4 \\ 
\hline\newstrut
10.134 & abbabbbcbAbbC &  6 & -13 &  -1 &  -3 & -18 &  7 &  -8 &  -8 &  -9 &  2 & -17 &  0 \\ 
\hline\newstrut
10.135 & abbcBAAABcc &  3 &  -1 &  3 &  -7 &  2 &  2 &  1 &  0 & -10 &  3 &  -7 &  1 \\ 
\hline\newstrut
10.136 & aabCbAbbCBB &  0 &  -1 &  -1 &  -1 &  4 &  3 &  2 &  8 &  -1 &  2 &  -3 &  -1 \\ 
\hline\newstrut
10.137 & aacBaBcADcBcBd & -2 &  2 &  -4 &  5 &  4 &  -3 &  2 &  -8 &  7 &  -6 &  4 &  0 \\ 
\hline
\end{tabular}}
\newpage

\cl{\begin{tabular}{|l|l|r|r|rr|rrr|rrrrr|} \hline\newstrut
Knot & Braid word  & $v_2$ & $v_3$ & $v_{4a}$ & $v_{4b}$ & $v_{5a}$ & $v_{5b}$ & $v_{5c}$ & $v_{6a}$ & $v_{6b}$ & $v_{6c}$ & $v_{6d}$ & $v_{6e}$ \\
\hline
\hline\newstrut
10.138 & aacBaBcdACbCbD & -3 &  2 &  -8 &  7 &  9 &  -2 &  5 & -11 &  10 & -13 &  0 &  -3 \\ 
\hline\newstrut
10.139 & aabbbaabab &  9 & -25 & -13 &  9 & -15 &  5 &  -7 & -23 & -11 &  1 &  2 &  -7 \\ 
\hline\newstrut
10.140 & aBaCCCbcccb &  2 &  -4 &  2 &  -2 &  0 &  6 &  0 &  4 & -10 &  0 &  -7 &  5 \\ 
\hline\newstrut
10.141 & AAAbaaBaab & -1 &  1 &  -5 &  -2 &  -2 &  -2 &  -3 &  0 &  2 &  -1 &  -1 &  -7 \\ 
\hline\newstrut
10.142 & Acccbcccaab &  8 & -21 & -11 &  1 & -14 &  9 & -12 & -20 & -11 & -10 &  -4 &  4 \\ 
\hline\newstrut
10.143 & AbbAAbaabb &  3 &  -5 &  3 &  -3 &  -2 &  6 &  0 &  2 & -12 &  0 & -10 &  4 \\ 
\hline\newstrut
10.144 & aacbbAbCCaB & -2 &  2 &  -8 &  0 &  -2 &  -4 &  -4 &  -8 &  4 &  -4 &  1 &  -7 \\ 
\hline\newstrut
10.145 & aabACbacbbc &  5 & -12 &  -4 &  -4 & -10 &  10 & -10 &  -8 & -14 & -12 & -12 &  10 \\ 
\hline\newstrut
10.146 & aaBCCbbaBcB &  0 &  0 &  2 &  3 &  8 &  2 &  5 & -12 &  0 &  -5 &  4 &  4 \\ 
\hline\newstrut
10.147 & ABBccbbaaBACb & -1 &  0 &  -4 &  0 &  3 &  1 &  1 &  1 &  3 &  -3 &  -2 &  -5 \\ 
\hline\newstrut
10.148 & aabbbAAbAb &  4 &  -7 &  1 &  -6 &  -7 &  7 &  -5 &  -3 & -11 &  -3 & -12 &  3 \\ 
\hline\newstrut
10.149 & AAbbbaaabb &  2 &  -2 &  0 &  -9 & -10 &  2 & -10 &  -4 &  -6 &  2 &  -1 &  0 \\ 
\hline\newstrut
10.150 & aabbcbABccB &  1 &  -1 &  -1 &  -6 &  -2 &  2 &  -4 &  4 &  -4 &  0 &  -4 &  -4 \\ 
\hline\newstrut
10.151 & aaCbbacBBCb &  3 &  -4 &  2 &  -6 &  -7 &  4 &  -5 &  -3 &  -8 &  -1 &  -8 &  1 \\ 
\hline\newstrut
10.152 & aabbbaaabb &  7 & -15 &  -1 &  -3 & -21 &  7 &  -8 &  1 &  -9 &  8 & -23 &  -5 \\ 
\hline\newstrut
10.153 & aabCCBBaCbb &  4 &  -1 &  3 &  -9 &  -5 &  0 &  -3 &  -1 & -12 &  3 & -10 &  2 \\ 
\hline\newstrut
10.154 & aBccbbaacbb &  5 &  -9 &  1 &  -8 & -20 &  6 & -12 &  -8 &  -8 &  4 & -12 &  3 \\ 
\hline\newstrut
10.155 & aaabAAbAAb & -2 &  2 &  -8 &  0 &  6 &  -2 &  1 &  -4 &  2 &  -5 &  -1 &  -8 \\ 
\hline\newstrut
10.156 & AbcbbcAAbCb &  1 &  -1 &  3 &  0 &  -9 &  -1 &  -4 & -13 &  -1 &  -2 &  2 &  2 \\ 
\hline\newstrut
10.157 & AbbaaBaabb &  4 &  -8 &  -2 &  -8 & -18 &  6 & -15 & -16 &  -6 &  -3 &  -5 &  8 \\ 
\hline\newstrut
10.158 & aabcABBcBBc & -3 &  1 & -11 &  4 &  8 &  0 &  4 & -10 &  4 &  -8 &  -2 &  -4 \\ 
\hline\newstrut
10.159 & aBaBabbbaB &  2 &  -3 &  3 &  -2 &  -8 &  2 &  -3 &  -8 &  -4 &  -1 &  -4 &  1 \\ 
\hline\newstrut
10.160 & abccaaBaBaC &  3 &  -6 &  -2 &  -8 &  -8 &  7 &  -9 &  4 &  -7 &  -3 & -11 &  2 \\ 
\hline\newstrut
10.161 & Abaabbbaab &  7 & -18 &  -8 &  0 & -15 &  9 & -12 & -25 & -13 & -12 &  -5 &  9 \\ 
\hline\newstrut
10.163 & AAcBccabbcB & -3 &  4 & -12 &  0 &  10 &  -5 &  1 &  -6 &  -1 &  -7 &  -1 & -12 \\ 
\hline\newstrut
10.164 & aBcAbaaBaBc &  1 &  -2 &  2 &  0 &  -8 &  1 &  -4 & -12 &  -1 &  -4 &  0 &  2 \\ 
\hline\newstrut
10.165 & AbCCaBaBacb &  1 &  0 &  4 &  1 &  8 &  2 &  5 & -10 &  -6 &  -1 &  3 &  6 \\ 
\hline\newstrut
10.166 & aabCbAbccAb &  2 &  -3 &  -1 &  -8 & -14 &  2 & -12 & -14 &  0 &  -4 &  -2 &  -3 \\ 
\hline
\end{tabular}}\normalsize

\end{document}

%% file: gtmonout.tex
\def\ifplaintex{\expandafter\ifx\csname documentclass\endcsname\relax}

\ifplaintex 
\else
\fi

\def\gt{{\mathsurround=0pt\it $\cal G\mskip-2mu$eometry \&\ 
$\cal T\!\!$opology}}        

\def\gtm{{\mathsurround=0pt\it $\cal G\mskip-2mu$eometry \&\ 
$\cal T\!\!$opology $\cal M\mskip-1mu$onographs}}    

\def\gtp{{\mathsurround=0pt\it $\cal G\mskip-2mu$eometry \&\ 
$\cal T\!\!$opology $\cal P\!$ublications}}  

\def\recd{{\small Received:\qua\receiveddate\ifx\reviseddate\relax
\else\qquad Revised:\qua\reviseddate\fi\par}} 

\def\volumenumber#1{\def\thevolumenumber{#1}}
\def\volumeyear#1{\def\thevolumeyear{#1}}
\def\volumename#1{\def\thevolumename{#1}}
\def\papernumber#1{\def\thepapernumber{#1}}
\def\pagenumbers#1#2{\def\startpage{#1}\def\finishpage{#2}}
\def\published#1{\def\publishdate{#1}}
\def\received#1{\def\receiveddate{#1}}
\def\revised#1{\def\reviseddate{#1}}
\def\accepted#1{\def\accepteddate{#1}}

\def\shorttitle#1{\def\theshorttitle{#1}}

\let\\\par
\let\thevolumenumber\relax\let\thepapernumber\relax
\let\thevolumeyear\relax\let\startpage\relax
\let\finishpage\relax\let\publishdate\relax\let\receiveddate\relax
\let\reviseddate\relax\let\accepteddate\relax\let\theasciititle\relax
\let\theasciiauthors\relax
\let\theasciiabstract\relax

\let\theerratum\relax\let\theasciiemail\relax
\let\theshortauthors\relax\let\theshorttitle\relax

\def\startpage{1}\def\finishpage{15}\def\thepapernumber{77}

\volumenumber{2}
\volumename{Proceedings of the Kirbyfest}
\volumeyear{1999}

\long\def\maketitlep{   

\count0=\startpage

\gtm\nl        
{\small Volume \thevolumenumber: \thevolumename\nl 
\ifx\theerratum\relax\else Erratum \erratumnumber\nl\fi
Pages \startpage--\finishpage\nl}

\vglue 0.1truein   

{\parskip=0pt\leftskip 0pt plus 1fil\def\\{\par\smallskip}{\ifplaintex\large
\else\Large\fi\bf\thetitle}\par\medskip}   
\vglue 0.05truein 

{\parskip=0pt\leftskip 0pt plus 1fil\def\\{\par}{\sc\theauthors}
\par\medskip}%
 
\vglue 0.03truein 

{\small\leftskip 25pt\rightskip 25pt{\bf Abstract}\stdspace\theabstract

{\bf AMS Classification}\stdspace\theprimaryclass
\ifx\thesecondaryclass\relax\else; \thesecondaryclass\fi\par
{\bf Keywords}\stdspace \thekeywords\par}\vglue 7pt

}   

\font\phead=cmsl9 scaled 950
\font=cmsl9 scaled 1050
\font\pnum=cmbx10 scaled 913
\font\lnum=cmbx10 
\font\pfoot=cmsl9 scaled 950
\font=cmsl9 scaled 1050
\ifplaintex
\headline{\vbox to 0pt{\vskip -4.5mm\line{\small\phead\ifnum
\count0=\startpage ISSN 1464-8997 (on line)
1464-8989 (printed) \hfill {\pnum\folio}\else\ifodd\count0\def\\{ }%
\ifx\theshorttitle\relax\thetitle\else\theshorttitle\fi\hfill{\pnum\folio}
\else\def\\{ and }{\pnum\folio}\hfill\ifx\theshortauthors\relax\theauthors
\else\theshortauthors\fi\fi\fi}\vss}}
\footline{\vbox to 0pt{\vglue 0mm\line{\small\pfoot\ifnum\count0=\startpage
Published \publishdate:\qua\copyright\ \gtp\hfill\else
\gtm, Volume \thevolumenumber\ (\thevolumeyear)\hfill\fi}\vss
}}
\else
\makeatletter
\def\@oddhead{{\small\lhead\ifnum\count0=\startpage ISSN 1464-8997 (on line)
1464-8989 (printed) \hfill {\lnum\number\count0}\else\ifodd\count0
\def\\{ }\ifx\theshorttitle\relax \thetitle \else\theshorttitle\fi\hfill
{\lnum\number\count0}\else\def\\{ and }{\lnum\number\count0}
\hfill\ifx\theshortauthors\relax 
\theauthors\else\theshortauthors\fi\fi\fi}}\def\@evenhead{@oddhead}
\def\@oddfoot{\small\lfoot\ifnum\count0=\startpage Published \publishdate:\qua\copyright\ \gtp\hfill\else
\gtm, Volume \thevolumenumber\ (\thevolumeyear)\hfill\fi}
\def\@evenfoot{@oddfoot}
\makeatother
\fi

\let\maketitlepage\maketitlep
\let\makeshorttitle\maketitlepage
\let\maketitle\maketitlepage

\newwrite\gtoutfile
\long\gdef\makeheadfile{  
{\def\\{, }\def\s{ }
\immediate\openout\gtoutfile head.xxx
\immediate\write\gtoutfile{Proxy-for: \ifx\theasciiauthors\relax
\theauthors\else\theasciiauthors\fi\s<\ifx\theasciiemail\relax\theemail\else\theasciiemail\fi>}
\immediate\write\gtoutfile{\noexpand\\}
\immediate\write\gtoutfile{Authors: \ifx\theasciiauthors\relax
\theauthors\else\theasciiauthors\fi}
{\def\\{ }\immediate\write\gtoutfile{Title: \ifx\theasciititle\relax
\thetitle\else\theasciititle\fi}}
\immediate\write\gtoutfile{Subj-class: GT or SG, GR etc}
\immediate\write\gtoutfile{MSC-class: \theprimaryclass\ifx\thesecondaryclass\relax\else, \thesecondaryclass\fi}
\immediate\write\gtoutfile{Journal-ref: Geom. Topol. Monogr. \thevolumenumber\s
(\thevolumeyear) \startpage-\finishpage}
\immediate\write\gtoutfile{Comments: Published by Geometry and Topology Monographs at}
\immediate\write\gtoutfile{\s\s\s  http://www.maths.warwick.ac.uk/gt/GTMon\thevolumenumber/paper\thepapernumber.abs.html}
\immediate\write\gtoutfile{\noexpand\\}
\immediate\write\gtoutfile{}
\ifx\theasciiabstract\relax
\immediate\write\gtoutfile{\theabstract}\else
\immediate\write\gtoutfile{\theasciiabstract}\fi
\immediate\write\gtoutfile{}
\immediate\write\gtoutfile{\noexpand\\}
\immediate\write\gtoutfile{}
\immediate\closeout\gtoutfile}}  

\def\maketitlepage{\maketitlep\makeheadfile}
\let\makeshorttitle\maketitlepage
\let\maketitle\maketitlepage